\documentclass{amsart}
\usepackage{amsmath, amssymb}
\newcommand{\eeq}{\end{equation}}

\newcommand{\bZ} {\mathbb{Z}_p}

\newcommand{\var}{z}

\newcommand{\sRS}{{\mathcal{MP}_p[\var]}}

\newcommand{\NN}{\mathbb{N}}
\newcommand{\QQ}{\mathbb{Q}}

\newcommand{\ZZ}{\mathbb{Z}}

\newcommand{\bbb}{c}
\newcommand{\ccc}{c'}
\newcommand{\aaa}{b}
\newtheorem{thm}{Theorem}[section]

\newtheorem{defi}[thm]{Definition}

\newtheorem{lem}[thm]{Lemma}

\def\eqref#1{(\ref{#1})}
\def\slfrac#1#2{\hbox{\kern.1em %
 \raise.5ex\hbox{\the\scriptfont0 #1}\kern-.11em %
 /\kern-.15em\lower.25ex\hbox{\the\scriptfont0 #2}}}

\title[Zero sets of holonomic sequences]{On the set of zero coefficients of a function satisfying a linear differential equation}

\author{Jason P. Bell}
\address{ Department of Mathematics, Simon Fraser University\\ 8888 University Dr.\\
Burnaby, BC, V5A 1S6\\ Canada}
\email{jpb@math.sfu.ca}

\author{Stanley N. Burris}
\address{Department of Pure Mathematics\\ University of Waterloo\\
Waterloo, ON,  N2L 3G1\\
Canada}
\email{snburris@math.uwaterloo.ca}

\author{Karen Yeats}
\address{Department of Mathematics, Simon Fraser University\\
8888 University Dr.\\ Burnaby, BC, V5A 1S6\\
Canada}
\email{karen.yeats@math.sfu.ca}
\thanks{The authors thank NSERC for supporting this project.}

\begin{document}

\begin{abstract} Let $K$ be a field of characteristic zero and suppose that 
$f:\mathbb{N}\to K$ satisfies a recurrence of the form
$$
f(n)\ =\ \sum_{i=1}^d P_i(n) f(n-i),
$$ 
for $n$ sufficiently large, where $P_1(\var),\ldots ,P_d(\var)$ are polynomials in $K[\var]$.  
Given that $P_d(\var)$ is a nonzero constant polynomial, 
we show that the set of $n\in \NN$ for which $f(n)=0$ is 
a union of finitely many arithmetic progressions and a finite set.  
This generalizes the Skolem-Mahler-Lech theorem, which assumes that
$f(n)$ satisfies a linear recurrence.  
We discuss examples and connections to the set of zero coefficients of a power series satisfying a homogeneous linear differential equation with rational function coefficients.
\end{abstract}
\maketitle

%\setlength{\baselineskip}{1.2\baselineskip}

%
%*****************************************************
%
% Section 1: Introduction
%
%*****************************************************
%
\section{Introduction}

The Skolem-Mahler-Lech theorem is a well-known result that describes the set of solutions 
to an equation $f(n)=0$, where $f(n)$ is a sequence given by a linear recurrence. 
Throughout this paper, we take $\NN := \{0,1,2,\ldots\}$.
%********************************************
% Theorem 1.1
%********************************************
\begin{thm} (Skolem-Mahler-Lech) \label{thm: SML}
Let $K$ be a field of characteristic zero and let $f:\mathbb{N}\to K$ be a sequence 
satisfying a linear recurrence over $K$; that is, a recurrence of the form
$f(n)=\sum_{j=1}^d a_jf(n-j)$, for $n\ge d$, with $a_1,\ldots ,a_d\in K$.  
Then the set $\big\{n \in \NN : f(n)=0\big\}$ is a union of finitely many arithmetic \label{thm: SMLlr}
progressions and a finite set. 
\end{thm}
This theorem was first proved for linear recurrences over the rational numbers by Skolem \cite{Sk34}, 
and it was then proved for linear recurrences over the algebraic numbers by Mahler \cite{Mah35}.  
The version above was proved first by Lech \cite{Le53} and later by Mahler \cite{Mah56, Mah57}.  
Further background about linear recurrences and the Skolem-Mahler-Lech theorem can be found 
in the book by Everest, van der Poorten, Shparlinski, and Ward \cite{EVSW03}.
%\vskip 2mm
There are many different proofs and extensions of the Skolem-Mahler-Lech theorem in the
literature \cite{Be06, Be06b, BGT, Bez89, Han86, vdP89, vdPT75}---all known proofs of 
the Skolem-Mahler-Lech theorem use $p$-adic methods in some way, 
albeit sometimes in a disguised manner.  

Recall that if $K$ is a field, then $f:\mathbb{N}\to K$ satisfies a linear recurrence over $K$ if and only if 
$$
\sum_{n=0}^{\infty} f(n)\var^n\in K[[\var]]
$$ 
is the power series expansion of a rational function in $K(\var)$.    
Stanley \cite[Theorem 2.1]{Stan} shows that
$$
\{{\rm Rational ~power ~series}\}~\subseteq ~  \{{\rm Algebraic ~power ~series}\}~\subseteq~
  \{{\rm Differentiably~finite~power~ series}\}.
$$
Algebraic power series comprise the power series $F(\var)\in K[[\var]]$ which satisfy a nonzero 
polynomial equation $P\big(\var ,F(\var)\big)=0$ with $P(\var ,y)\in K[\var ,y]$.  
The set of \emph{differentiably finite} power series $F(\var)$ are those that satisfy a non-trivial 
homogeneous linear differential equation with rational function coefficients:
\begin{equation}
\sum_{i=0}^d Q_i(\var) \frac{d^i}{d\var ^i}F(\var) \ = \ 0.
\end{equation}
\begin{defi}{\em
Let $K$ be a field and let $f:\mathbb{N}\to K$ be a $K$-valued sequence.  
We say that $f(n)$ is \emph{holonomic} if the power series 
$$
\sum_{n=0}^{\infty} f(n)\var ^n\in K[[\var]]
$$ 
is differentiably finite.}
\end{defi}
It is straightforward to show that a $K$-valued sequence $f(n)$ is holonomic if 
it satisfies a \emph{polynomial-linear recurrence}; that is, there exist a natural number $d$ and polynomials 
$P_0(\var),\ldots ,P_d(\var)\in K[\var]$, not all zero, such that
\begin{equation}
\sum_{i=0}^d P_i(n)f(n-i) \ = \ 0,
\end{equation}
for all sufficiently large $n$ (see Stanley \cite[Theorem 1.5]{Stan}).  
It will be convenient to use this characterization of holonomic sequences throughout this paper.

Differentiably finite power series and holonomic sequences have been extensively studied by several
 authors \cite{CS, CMS, Garo, Ges, Stan, WZ}, and many important sequences arising naturally in
  combinatorics, number theory, and algebra have been shown to be holonomic 
  \cite{BM, Ges, LVO, vdPAp, Stan}.  
A number of significant results, including the following, have employed holonomic methods: 
   (i) Ap\'ery's proof of the irrationality of $\zeta(3)$ used a sequence satisfying a polynomial-linear
    recurrence  (see van der Poorten's survey \cite{vdPAp});
    (ii) Wilf and Zeilberger \cite{WZ} showed how the theory of holonomic sequences could be 
    applied to the ``automatic'' proving of combinatorial identities; and 
    (iii) Bousquet-M\'elou \cite{BM} applied the theory to enumeration of walks in the quarter plane.

In light of the well-behaved nature of the zero set of a sequence satisfying a linear recurrence 
over a field of characteristic zero, 
it is natural to ask whether a similar result holds for holonomic sequences.  
Indeed, Rubel \cite[Problem 16]{Rub} asked: {\em does the conclusion of the Skolem-Mahler-Lech 
theorem hold for all holonomic sequences?}  The strongest result in this direction is due to 
B\'ezivin \cite{Bez89} and Methfezzel \cite{Meth}: the set of zeros of a 
holonomic sequence can be expressed as a union of finitely many arithmetic progressions and a 
set of zero density.  (B\'ezivin assumed that $0$ and $\infty$ are not irregular 
singular points of the corresponding linear differential equation.)  
Laohakosol \cite{Lao} and B\'ezivin and Laohakosol \cite{BL} showed that the answer to
Rubel's question is `yes',  provided additional technical conditions on the associated differential equation hold.  
(There is an error in the first paper of Laohakosol 
which is repaired in the subsequent paper with B\'ezivin.)

We are able to give an affirmative answer to Rubel's question in the case that the recurrence has 
a reciprocal property (see Stanley \cite[\S 3]{Stan}), namely the recurrence can be `run backwards'
in a well-behaved way.  

%********************************************
% Theorem 1.3
%********************************************
\begin{thm} \label{thm: main} 
Let $K$ be a field of characteristic zero, 
let $d$ be a positive integer, and let $P_1(\var),\ldots ,P_d(\var)\in K[\var]$ be polynomials 
with $P_d(\var)$ a nonzero constant.   
Suppose that $f:\mathbb{N}\to K$ is a sequence satisfying the polynomial-linear recurrence
\begin{equation} \label{poly linear}
f(n)\ =\ \sum_{i=1}^d P_i(n)f(n-i)
\end{equation} 
for all $n$ sufficiently large.  
Then $\big\{n\in \mathbb{N}~:~ f(n)=0\big\}$ is a union of finitely many arithmetic progressions 
and a finite set.
\end{thm}
This theorem generalizes the Skolem-Mahler-Lech theorem, which asserts that 
the conclusion holds provided $P_1(\var),\ldots ,P_d(\var)$ are constant polynomials.
%\footnote{
%FOOTNOTE 
We note that the reciprocal property of Stanley mentioned above is given by the assumption that $P_d(\var)$ is a nonzero 
constant polynomial---this allows the recurrence to be rewritten as
$$
f(n)\ =\ \sum_{i=1}^{d-1} \lambda P_{d-i}(n+d)f(n+i) -\lambda f(n+d),
$$ 
where $\lambda = -1/P_d(\var)$, a nonzero element of $K$.
Thus the definition of $f(n)$ extends to all $n\in\mathbb{Z}$. 
%}
%END FOOTNOTE 

The proof of Theorem \ref{thm: main} uses $p$-adic methods and consists of three steps.  
The first step uses an argument of Lech to show that one can assume 
the base field $K$ is $\mathbb{Q}_p$, for some prime $p$, and the nonzero 
coefficients of the polynomials $P_i(\var)$ in \eqref{poly linear} are units in the $p$-adic integers. 
The next step involves showing that there exists a positive integer $\aaa$ such that 
for $\bbb\in \{0,\ldots ,\aaa-1\}$, the sequence $f(\aaa n+\bbb)$ can be embedded in an `analytic arc', 
in the sense that there exists a $p$-adic analytic map $g:\bZ\to \bZ$ such that $g(n)=f(\aaa n+\bbb)$.  
We then use a theorem of Strassman (see \S 2) to show that 
$\big\{n \in \NN : f(\aaa n+\bbb) = 0\big\}$ is either finite or $\NN$.
The difficult step is the second step, in which one must find a $p$-adic analytic map that agrees with a subsequence on the natural numbers.
%\footnote{
%FOOTNOTE
 In the case that $f(n)$ satisfies a linear recurrence over a field $K$, 
it is well-known (see, for example, Everest, van der Poorten, Shparlinski, and Ward \cite[Section 1.1.6]{EVSW03}) that there exist constants $\gamma_{i,j}\in \overline{K}$ and 
$\alpha_1,\ldots ,\alpha_m\in \overline{K}$ such that, for some positive integer $N$ and all $n$ sufficiently large, we have
$$
f(n)\ =\ \sum_{i=1}^m \sum_{j=0}^N \gamma_{i,j} n^j \alpha_i^n.
$$  
Lech used this fact to show that if $K=\mathbb{C}_p$ and $|\alpha_i-1|_p<p^{-1/(p-1)}$ 
then the map
$$
g(z)\ :=\ \sum_{i=1}^m \sum_{j=0}^N \gamma_{i,j} z ^j \alpha_i^z 
$$ 
will be an analytic map that converges absolutely on the closed unit ball of $\mathbb{C}_p$.  
Such an expression for $f(n)$ is not available when one works with holonomic sequences, and 
additional work is needed to obtain the embedding into a $p$-adic analytic arc.  
%}
%END FOOTNOTE

\section{Preliminaries}
This section has some of the $p$-adic results necessary to prove Theorem \ref{thm: main}.  
We start with an embedding theorem due to Lech \cite{Le53}---this result can be regarded 
as a $p$-adic analogue of the Lefschetz principle.
%********************************************
% Lemma
%********************************************
\begin{lem}~\label{lem: embed}  
Let $K_0$ be a finitely generated extension of $\mathbb{Q}$ and 
let $\mathcal{S}$ be a finite subset of $K_0$.  
Then there exist infinitely many primes $p$ such that 
$K_0$ embeds in $\mathbb{Q}_p$; moreover, for all but
finitely many of these primes, 
every nonzero element of $\mathcal{S}$ 
is sent to a unit in $\mathbb{Z}_p$.
\end{lem}
\begin{proof} 
The first-named author \cite[Lemma 3.1]{Be06} used the  Chebotarev density theorem 
(see Lang \cite[Theorem 10, p. 169]{La94}) 
to show that the set of primes $p$ having the required property has positive density.
\end{proof}
The other main result from $p$-adic analysis required is Strassman's theorem \cite{St28}, 
which  asserts that if a power series $f(\var)\in \mathbb{Q}_p[[\var]]$  converges in the closed 
$p$-adic unit disc
$$
\overline{B_{\QQ_p}(0; 1)} \ :=\  \big\{ \alpha \in \QQ_p:~ |\alpha|_p \le 1 \big\}\quad ( =\ \ZZ_p),
$$
and has infinitely many zeros in this disk, then it is  identically zero.  
A power series $f(\var)\in \mathbb{Q}_p[[\var]]$ that converges in the closed 
$p$-adic unit disc will be called a \emph{rigid $p$-adic power series}.

To prove Theorem \ref{thm: main}, we use the following ring $\sRS$ of polynomials. 
%********************************************
% Defn. 6.1/ 5.4
%********************************************
\begin{defi}~\label{de44}
{\em 
Given a prime $p$, let $\sRS$ be the subring of the polynomial ring
${\mathbb Q}_p[\var]$ consisting of the polynomials 
which map $\ZZ_p$ into itself.
}
\end{defi}
\smallskip
Clearly $\mathbb{Z}_p[\var] \subseteq \sRS \subseteq \mathbb{Q}_p[\var]$.
A theorem of Mahler \cite[p. 49-50]{Mah81} shows that for $p$ a prime,  
\begin{equation}
\sRS \ =  \ 
\Bigg\{ 
\sum_{i=0}^m \alpha_i {\var \choose i}~:~ m\ge 0~{\rm and}~\alpha_i\in \ZZ_p \Bigg\}.
\end{equation}

%%%%%%%%
%%%%%%
\section{Proof of Theorem \ref{thm: main}}

In this section we prove our main result.  We begin with a simple lemma.

\begin{lem} \label{lem: matrix}
Let $p$ be prime and suppose that $Q_1(\var), \ldots, Q_{d}(\var) \in \sRS$  
are such that each is congruent modulo $p\ \sRS$ to a polynomial in $\sRS$ 
of degree at most $N$.
Then there exist 
$H_1(\var),\ldots ,H_d(\var)\in \sRS$, each of degree at most $N+1$, with 
$H_1(0)=\cdots =H_d(0)=0$,
such that
\begin{equation} \label{eq: mateq} 
\left[
\begin{array}{c} 
H_{1}(\var +1) \\ \vdots \\ H_{d}(\var +1)
\end{array}
\right]
\ \equiv \ 
\left[
\begin{array}{c} 
H_{1}(\var) \\ \vdots \\ H_{d}(\var)
\end{array}
\right] 
-
\left[
\begin{array}{c} 
Q_{1}(\var) \\ \vdots \\ Q_{d}(\var)
\end{array}
\right]
\ \big({\rm mod}~p\sRS \big).
\end{equation}
%
%where $H_1(0)=\cdots =H_d(0)=0$ and $H_1,\ldots ,H_d$ have degree at most $N+1$.
\end{lem}
\begin{proof} It is no loss of generality to assume that 
$Q_1(\var), \ldots, Q_{d}(\var)$ each have degree at most $N$. 
Thus 
\[ 
Q_{i}(\var) \ = \ \sum_{k=0}^{N} \alpha_{i,k} {\var\choose k}
\]
with each $\alpha_{i,k} \in \ZZ_p$.
Let 
\begin{equation} \label{eq: hi0}
H_{i}(\var) \ := \ - \sum_{k=0}^N \alpha_{i,k} {\var \choose k+1}.
\end{equation}
Then $H_{i}(\var) \in \sRS$, and it is of degree at most $N+1$.
It is easy to check that this gives a solution 
to equation (\ref{eq: mateq}), using the identity
\[
 {\var +1\choose k+1}-{\var  \choose k+1} \ =\  {\var  \choose k}.
\]
Furthermore, $H_i(0)=0$ for $1\le i\le d$.
\end{proof}

To create analytic maps for the modified version of the $p$-adic analytic arc lemma, 
we will use the following lemma about a subalgebra of $\sRS$. 
%The ring $\sRSp$ is the set of
%polynomials $P(p\var )$, for $P(\var) \in \sRS$.

%\
% USED to bound degrees of the h's in the analytic arc lemma
%
%LEMMA
\begin{lem} Given a prime $p$ and a positive integer $m$, let
\[ 
S_m \ := \ \Bigg\{ \gamma + \sum_{i=1}^m p^i H_i(\var)~\Big|~\gamma\in \ZZ_p, H_i(\var)\in \sRS,
 {\rm deg}\big(H_i(\var)\big)\le 2i-1\Bigg\}
 \]
and
\[ 
T_m \ := \ S_m+\Bigg\{ \gamma + \sum_{i=1}^n p^i H_i(\var)~\Big|~n\ge 1, \gamma\in \ZZ_p, 
H_i(\var)\in \sRS, {\rm deg}\big(H_i(\var)\big)\le 2i-2\Bigg\}.
\]
Then  $T_m$ is a $\ZZ_p$-subalgebra of $\sRS$.
%,  $\sRSp \subseteq T_m$,
Moreover, if $\alpha\in p\ZZ_p$ and $\beta\in \ZZ_p$, then for any polynomial $P(\var)\in \ZZ_p[\var]$, 
one has $P(\beta+\alpha \var ) \in T_m$, for $m\ge 1$.   In particular, $\ZZ_p[p\var]\subseteq T_m$ for $m\ge 1$.
\label{lem: SN}
\end{lem}
\begin{proof} Since $T_m$ is closed under addition and multiplication by $p$-adic integers,
to show $T_m$ is a $\ZZ_p$-algebra, 
it is sufficient to show that $T_m$ is closed under multiplication. Since $T_m$ is contained in the $\ZZ_p$-span of the constant function $1$ and elements of the form $p^i H(\var)$ for some $i\ge 1$ and $H(\var)\in \sRS$ a polynomial of degree at most $2i-1$, it suffices 
to check that if $i$ and $j$ are positive integers and $H(\var),G(\var)\in \sRS$ are polynomials of degrees at most $2i-1$ and $2j-1$ respectively, then the product $p^i H(\var) \cdot p^j G(\var)$ lies in $T_m$.  Note that this follows since $H(\var)G(\var)\in \sRS$ has degree at most $2(i+j)-2$ and so the product is of the form $p^{i+j}L(\var)$ for some $L(\var)\in \sRS$ of degree at most $2(i+j)-2$.  

Next suppose $m \ge 1, \alpha\in p\ZZ_p$, $\beta\in \ZZ_p$, and $P(\var)\in \ZZ_p[\var]$. 
Then $\alpha$ in $p\ZZ_p$ implies
$\beta+\alpha \var  \in T_m$; and since $T_m$ is a $\ZZ_p$-algebra, 
it follows that $P(\beta+\alpha \var ) \in T_m$, and so, in particular, $\ZZ_p[pz]\subseteq T_m$.
\end{proof}

%************************
%
% analytic arc lemma 
%*************************
Given a positive integer $d$ and a ring $R$ with identity element, we let $\mathsf{M}_d(R)$ denote the set of $d\times d$ matrices
with entries from $R$.  When working with a matrix ring $\mathsf{M}_d(R)$, we take ${\bf I}$ to be the $d\times d$ identity matrix.
%LEMMA
\begin{lem} (Analytic arc lemma) 
Let $d$ be a natural number, let $p\ge 5$ be prime, let 
${\bf v}=[v_1\,\cdots \, v_d]^{ T} \in \ZZ_p^d$, and let 
$\mathbf{A}(\var):=\big(a_{ij}(\var)\big) \in \mathsf{M}_d\big(\ZZ_p[p\var ]\big)$ 
with 
$\mathbf{A}(0) \equiv \mathbf{I}\ \left(\bmod\, \ p\mathsf{M}_d\big(\ZZ_p\right)$.
%$a_{ij}(\var) \in \sRSp$.
Then there exist rigid $p$-adic power series $f_1(\var),\ldots ,f_d(\var)$
%:\mathbb{Z}_p\to \mathbb{Z}_p$ 
such that $f_i(0)=v_i$ for $i\in \{1,\ldots ,d\}$ and
\[
\left[ \begin{array}{c} f_1(\var +1) \\ \vdots \\ f_d(\var +1)\end{array}\right] \ =\ \mathbf{A}(\var)\cdot
\left[ \begin{array}{c} f_1(\var) \\ \vdots \\ f_d(\var)\end{array}\right].
\] 
%for every $x\in\mathbb{Z}_p$. 
\label{lem: aal}
\end{lem}

\begin{proof} Let $S_m$ and $T_m$ be as in the statement of Lemma \ref{lem: SN}.
The desired tuple $\big(f_1(\var), \ldots, f_{d}(\var)\big)$ of power series will be successively 
approximated by tuples of polynomials  $\big(G_{1, j}(\var), \ldots, G_{d,j}(\var)\big)$, that is,
$f_i(\var) - G_{i,j}(\var) \in p^{j} \sRS$, for $1 \le i \le d$ and $j\ge 0$.
 Let 
$$
G_{i,0}(\var) \ :=\  v_i   ~~~~~\mbox{for}~~ 1\le i \le d.
$$
It will be proved, by induction on $m$, that one can recursively find $H_{i,m}(\var) \in \sRS$, for $1\le i\le d$,
such that by setting
\begin{equation} \label{def Gim}
G_{i,m}(\var) \ :=\ v_i+\sum_{k=1}^{m} p^{k}H_{i,k}(\var),
\end{equation}
one has the following three conditions holding:
\begin{enumerate}
\item[(i)] $H_{i,m}(0)=0$ for $1\le i\le d$;
\item[(ii)] $H_{i,m}(\var)\in S_m$, for  $1\le i\le d$;
\item[(iii)] and 
$$ 
G_{i,m}(\var +1)\ - \sum_{j=1}^d a_{i,j}(\var)G_{j,m}(\var)\ \in\ p^{m+1} \sRS.
$$
\end{enumerate}

 The base case of the induction is $m=0$.  In this case, conditions (i) and (ii) are
 vacuous, and (iii) holds since $\mathbf{A}(\var)$ is congruent modulo $p$ to the identity matrix.

  Let $m \ge 0$ and assume that the $H_{i,k}(\var)$ have been found, for $1\le i\le d$ 
and $0 \le k \le m $, such that conditions (i)--(iii) hold for $0 \le k \le m$.  
The method is to find polynomials $H_{i, m+1}(\var) \in \sRS$
such that, with $G_{i,m+1}(\var)$ defined as in \eqref{def Gim},
conditions (i)--(iii) hold. 

By (iii), there are polynomials $Q_{i,m}(\var)\in \sRS$,
for $1\le i\le d$, such that
\begin{equation} \label{eq: gH}
p^{m+1} Q_{i,m}(\var)\ = \ G_{i,m}(\var +1)-\sum_{j=1}^d a_{i,j}(\var)G_{j,m}(\var).
%G_{i,m}(x+1)-\sum_{j=1}^d a_{i,j}(\var)G_{j,m}(\var) =p^{m+1}Q_{i,m}(\var).
\end{equation}
%
%we see that 
Then Definition \eqref{def Gim} and condition (ii) show that
$G_{1,m}(\var),\ldots ,G_{d,m}(\var)$ as well as $$G_{1,m}(\var +1),\ldots ,G_{d,m}(\var +1)$$ 
are in $S_{m}$.
Thus, by Lemma \ref{lem: SN} and the fact that the $a_{ij}(\var)$ are in
$\ZZ_p[pz]$, $ p^{m+1} Q_{i,m}(\var)$ is in the $\ZZ_p$-algebra $T_{m}$.  
It follows that 
$$
p^{m+1} Q_{i,m}(\var)\ =\  \gamma_{i,m}+ \sum_{k=1}^n p^k Q_{i,m,k}(\var)
$$ 
for some $\gamma_{i,m}\in  \mathbb{Z}_p$, and for some polynomials $Q_{i,m,k}(\var)\in \sRS$ 
such that ${\rm deg}\big(Q_{i,m,k}(\var)\big)\le 2k-1$ for $k\le m$ and 
${\rm deg}\big(Q_{i,m,k}(\var)\big)\le 2k-2$ for $k> m$.  

Consequently,
$p^{m+1}Q_{i,m}(\var)$ is equivalent modulo $p^{m+2}\sRS$ to the polynomial
$$
 \gamma_{i,m}+ \sum_{k=1}^{m+1} p^k Q_{i,m,k}(\var),
 $$
a polynomial in $\sRS$ of degree at most $2m$.   
Hence $Q_{i,m}(\var)$ is congruent modulo $p\sRS$ to a polynomial in $\sRS$ 
of degree at most $2m$.

Note that the definition \eqref{def Gim} of $G_{i,m+1}(\var)$ can be replaced by
\[
G_{i, m+1}(\var) \ :=\ G_{i,m}(\var) + p^{m+1} H_{i, m+1}(\var).
\]
To satisfy property (iii), it  is sufficient to find 
$H_{i, m+1}(\var) \in \sRS$, $ 1 \le i \le d$, such that 
\[
G_{i,m}(\var +1)+p^{m+1} H_{i,m+1}(\var +1)  - 
\sum_{j=1}^d a_{i,j}(\var)\big( G_{j,m}(\var) +p^{m+1} H_{j,m+1}(\var)\big) 
\]
is in  $p^{m+2} \sRS$, for $ i\in \{1,\ldots , d\}$.
This expression is congruent modulo $p^{m+2} \sRS$ to
\[
p^{m+1}Q_{i,m}(\var) + p^{m+1} H_{i,m+1}(\var +1)
 - p^{m+1}\sum_{j=1}^d a_{i,j}(\var) H_{j,m+1}(\var).
 \]
However, since each $a_{i,j}(\var)\in \delta_{i,j}+p \sRS$, we see that this simplifies as
\[
p^{m+1}Q_{i,m}(\var) + p^{m+1} H_{i,m+1}(\var +1)
 - p^{m+1} H_{i,m+1}(\var)
 \]
modulo $p^{m+2}\sRS$. 
It therefore suffices  to solve the system
\begin{equation} \label{eq: Q}
Q_{i,m}(\var) + H_{i,m+1}(\var +1)
 - H_{i,m+1}(\var) \ \equiv\  0 \ \left(\bmod \, p \sRS\right),
\end{equation}
for $i\in \{1,\ldots ,d\}$, where each $Q_{i,m}$ is congruent modulo $p\sRS$ to a polynomial of degree at most $2m$.

The hypotheses of the statement of Lemma \ref{lem: matrix} are satisfied by the 
system \eqref{eq: Q}, so
we conclude that there exists a  solution 
$\big[ H_{1,m+1}(\var), \ldots, H_{d, m+1}(\var) \big] \in {\sRS}^d$ 
with $H_{i,m+1}(0)=0$ for $1\le i\le d$ 
and with
$H_{i,m+1}(\var)$ of degree at most $2(m+1)-1$.
Thus conditions (i)--(iii) are satisfied, completing the induction
step.

We set
\[
f_i(\var) \ :=\  v_i + \sum_{j=1}^{\infty} p^j H_{i, j}(\var).
\]
Then each $H_{i,j}(\var) \in \sRS$ is of degree at most $2j-1$ and hence
\[
H_{i,j}(\var)\ =\ \sum_{k=0}^{2j-1}  \gamma_{i,j,k} {\var  \choose k},
\]
with $\gamma_{i,j,k} \in \mathbb{Z}_p$.  (Let $\gamma_{i,j,k}=0$ for $k>2j-1$.) 
We find that
\begin{eqnarray}
f_i(\var) &=& v_i + \sum_{j=1}^{\infty} p^j 
\left( \sum_{k=0}^{2j-1} \gamma_{i,j,k} {\var  \choose k}\right) \nonumber\\
&=& v_i + \sum_{k=0}^{\infty} \beta_{i,k}  {\var  \choose k}, \label{eq475}
\end{eqnarray}
in which
\[
\beta_{i,k} \ :=\  \sum_{j=1}^{\infty} p^j \gamma_{i,j,k}
\]
is absolutely convergent $p$-adically, since each $\gamma_{i,j,k} \in \mathbb{Z}_p$.
To show that the series (\ref{eq475}) defines an analytic function on
 $\ZZ_p$, we must establish that $|\beta_{i,k}|_{p}/|k!|_p \to 0$
as $k \to \infty$  (see Robert \cite[Theorem 4.7, p. 354]{Rob}); that is, for any $j >0$ one has $\beta_{i,k}/k! \in p^j \bZ$
for all sufficiently large $k$.
To do this, we note that $\gamma_{i,j,k}=0$ if $j<(k+1)/2$.  Hence
\[
 \beta_{i,k} \ =\  \sum_{j\ge (k+1)/2} p^j \gamma_{i,j,k}.
 \]
It follows that $|\beta_{i,k}|_p\le p^{-(k+1)/2}$. 
Since $1/|k!|_p<p^{k/(p-1)}$, we see 
that 
$\beta_{i,k}/k!\rightarrow 0$ since $p>3$.  
Hence $f_1,\ldots ,f_d$ are rigid analytic maps on $\bZ$.
 
Finally, observe that the argument above showed that
\[ 
f_i(\var) \ \equiv\  G_{i,j}(\var)~\big(\mathrm{mod}~p^{j}\sRS\big).
\]
It then follows from property (iii) above that
\[
f_i(\var +1) \ \equiv\  \sum_{\ell=1}^d a_{i,\ell}(\var) f_{\ell}(\var)\  
\big(\mathrm{mod}~ p^j \sRS\big)
\]
for $i\in \{1,\ldots , d\}$.

Since this holds for all $j \ge 1$, we conclude
that
\[
\left[ \begin{array}{c} f_1(\var +1) \\ \vdots \\ f_d(\var +1)\end{array}\right] \ = 
\ \mathbf{A}(\var)\cdot \left[ \begin{array}{c} f_1(\var) \\ \vdots \\ f_d(\var)\end{array}\right].
\]
Finally, we have
\[ 
f_{i}(0)\ =\ v_i + \sum_{j=1}^{\infty} p^j H_{i,j}(0)\ =\ v_i,
\]
which concludes the proof.
\end{proof}
We are almost ready to prove our main result.  The one remaining thing we need is to 
show how the analytic arc lemma applies to our situation.  
We accomplish this with the following lemma.

%LEMMA
\begin{lem} Let $p$ be a prime number, let $d$ be a natural number, and 
let $P_1(\var),\ldots ,P_d(\var)\in \bZ[\var]$.   
Suppose that $f:\mathbb{N}\to \bZ$ is a sequence satisfying the polynomial-linear recurrence
$$
f(n)\ =\ \sum_{i=1}^d P_i(n)f(n-i)
$$ 
for $n\ge d$. 
Let
\[
\mathbf{B}(\var) \ := \ \left( \begin{array}{ccccc} 0 & 1 & 0 & \cdots & 0 \\
0 & 0 & 1&  \cdots &  0 \\
\vdots & \vdots & \ddots & \vdots & \vdots \\
0 & 0 & \cdots & 0 & 1 \\
P_d(\var -1+d) & P_{d-1}(\var -1+d) & \cdots &  P_2(\var -1+d) & P_1(\var -1+d)\end{array}\right), 
\]
and let ${\bf v}_0 := [f(0), \ldots , f(d-1)]^{T}$. 
Then 
%$\mathbf{B}(\var) \in \mathsf{M}_d\big(\ZZ_p[\var]\big)$, ${\bf v}_0^T \in {\ZZ_p}^d$, 
%and
\[
\left[ 
\begin{array}{c} f(n) \\ \vdots \\ f(n+d-1)\end{array}
\right] 
 \ =\ \mathbf{B}(n)\mathbf{B}(n-1)\cdots \mathbf{B}(1){\bf v}_0,
\] 
for every positive integer $n$. 
 Furthermore, if $P_d(\var)$ is a constant that is a unit in $\bZ$ then 
 the determinant of $\mathbf{B}(\var)$ is a constant that is a unit in $\bZ$.
\label{lem: mat1}
\end{lem}
\begin{proof}
For each positive integer $n$ let
\[
{\bf v}_n \ :=\ \left[ \begin{array}{c} f(n) \\ \vdots \\ f(n+d-1)\end{array}\right] .
\]
Then
$$
{\bf v}_n\ =\ \mathbf{B}(n){\bf v}_{n-1},
$$ 
for $n\ge 1$,
leads to
$$
{\bf v}_n\ =\ \mathbf{B}(n)\mathbf{B}(n-1)\cdots \mathbf{B}(1){\bf v}_0.
$$
The determinant of $\mathbf{B}(\var)$ is $P_d(\var -1+d)$; hence 
if $P_d(\var) = \alpha$, a unit in $\ZZ_p$, then $\det\big(\mathbf{B}(\var)\big)$ equals
that unit.
\end{proof}

\begin{proof}[Proof of Theorem \ref{thm: main}]

It is no loss of generality to assume that the recurrence in the statement of Theorem \ref{thm: main} holds for all $n\ge d$.  Let $\mathcal{S}$ denote the finite subset of $K$ consisting of the nonzero $f(n)$, 
where $0 \le n < d$, along with all nonzero coefficients of the polynomials 
$P_1(\var),\ldots ,P_d(\var)$.  
Then we let $K_0 = \mathbb{Q}(\mathcal{S})$, 
the subfield of $K$ generated by the elements of $\mathcal{S}$.    
By Lemma \ref{lem: embed}, there exists a prime $p\ge 5$ and a field embedding of $K_0$ into 
$\mathbb{Q}_p$ such that the image of $\mathcal{S}$ is contained in the units of $\mathbb{Z}_p$.
By identifying $K_0$ with its image in $\mathbb{Q}_p$, we see that there is no loss of generality 
in assuming that $P_1(\var),\ldots ,P_d(\var)\in \mathbb{Z}_p[\var]$, with the nonzero coefficients of the $P_i(\var)$ being
units of $\ZZ_p$; in particular, $P_d(\var)$ is a unit in $\ZZ_p$. 
Furthermore, after this identification, a simple induction shows that $f : \mathbb{N} \rightarrow \ZZ_p$.

Letting $\mathbf{B}(\var)$ and ${\bf v}_0$ be as defined in the statement of Lemma \ref{lem: mat1}, 
one has, for $n \ge d$,
\begin{equation} \label{cascade}
\left[ 
\begin{array}{c}
 f(n) \\ \vdots \\ f(n+d-1)
 \end{array}
 \right] 
\ =\ 
\mathbf{B}(n)\mathbf{B}(n-1)\cdots \mathbf{B}(1){\bf v}_0,
\end{equation}
and  the determinant of $\mathbf{B}(\var)$ is a unit in $\mathbb{Z}_p$.  
In addition, observe that $(\var +p)^i-\var ^i \in p\bZ[\var]$ for all nonnegative integers $i$ and 
hence, by linearity, we have that $Q(\var +p)-Q(\var)\in p\bZ[\var]$ for all $Q(\var)\in \bZ[\var]$.  
Thus we have $\mathbf{B}(\var +p)-\mathbf{B}(\var)\in p \mathsf{M}_d\big(\bZ[\var]\big)$, 
since the entries of $\mathbf{B}(\var)$ are in $\bZ[\var]$.

Let $\varphi: M_d(\bZ)\to M_d(\mathbb{Z}/p\mathbb{Z})$ be the canonical surjection.   For $n\in \ZZ$, $\mathbf{B}(n) \in \mathsf{M}_d(\ZZ_p)$, and thus 
$\varphi\big(\mathbf{B}(n)\big) \in \mathsf{M}_d(\ZZ/p\ZZ)$. 
Since $\mathsf{M}_d(\ZZ/p\ZZ)$ is a finite ring, there exist natural numbers $m_0$ and $m_1$ with $m_0<m_1$ such that
$$
\varphi\big(\mathbf{B}(pm_1)\cdots \mathbf{B}(1)\big) \ = 
\ \varphi\big(\mathbf{B}(pm_0)\cdots \mathbf{B}(1)\big).
$$
Since the determinant of $\mathbf{B}(\var)$ is a unit of $\ZZ_p$, the determinant of 
$\varphi\big(\mathbf{B}(i)\big)$ is nonzero, for
$i \in \ZZ$. Thus $\varphi\big(\mathbf{B}(i)\big)$ is in ${\rm GL}_d(\ZZ/p\ZZ)$ and hence
$\varphi\big(\mathbf{B}(pm_1)\cdots \mathbf{B}(pm_0+1)\big)$ is the identity.  
Furthermore, since $\varphi\big(\mathbf{B}(n+p)\big) =\varphi\big(\mathbf{B}(n)\big)$ for $n\in\ZZ$, 
it follows that for $i\in \ZZ$, $\varphi\big(\mathbf{B}(pm_1+i)\cdots \mathbf{B}(pm_0+1+i)\big)$ is 
similar to $\varphi\big(\mathbf{B}(pm_1)\cdots \mathbf{B}(pm_0+1)\big)$ and hence
\begin{equation}\label{eq: x}
\varphi\big(\mathbf{B}(pm_1+i)\cdots \mathbf{B}(pm_0+1+i)\big)\ =\ {\bf I}.
\end{equation}
Let
\begin{equation*}
 {\aaa} \ :=\  p(m_1-m_0).
\end{equation*}

Suppose that $\bbb\in \{0,\ldots ,{\aaa}-1\}$.
Then Equation (\ref{eq: x}) gives
\begin{equation}
\label{eq: y}
\varphi\big(\mathbf{B}(\bbb+{\aaa} (n+1))\cdots \mathbf{B}(\bbb+{\aaa} n+1)\big)\ =\ {\bf I},
\end{equation}
for  $n \in \ZZ$.
Define the matrix $\mathbf{A}(\var)\in \mathsf{M}_d\big(\bZ[\var]\big)$ by
\begin{equation*}\label{eq: xx}
\mathbf{A}(\var)\ :=\ \mathbf{B}(\bbb+{\aaa} \var )\mathbf{B}(\bbb+{\aaa} \var -1)\cdots 
\mathbf{B}(\bbb+{\aaa} \var -\aaa+1).
\end{equation*}
Since $p$ divides $a$, we see that for each integer $i$, the matrix-valued function $\mathbf{B}(i+\aaa z)\in \mathsf{M}_d\big(\bZ[p\var]\big)$ and hence $\mathbf{A}(\var)\in \mathsf{M}_d\big(\bZ[p \var]\big)$.  
Moreover, Equation (\ref{eq: y}) shows that for $n\in \mathbb{Z}$ we have $\varphi\big(\mathbf{A}(n)\big)={\bf I}$ and in particular $\mathbf{A}(0)\equiv {\bf I}~\left(\bmod\, p\mathsf{M}_d(\bZ)\right)$.  
%Since $\varphi$ is continuous and $\mathbb{Z}$ is dense in $\bZ$, 
%we have $\varphi\big(\mathbf{A}(\alpha)\big)={\bf I}$ for all $\alpha\in \mathbb{Z}_p$; 
%that is 
%$\mathbf{A}(\var)\in {\bf I}+p\cdot\mathsf{M}_d\big(\sRS\big)$.
Let 
\[
 {\bf v} \ :=\ \mathbf{B}(\bbb) \mathbf{B}(\bbb-1)\cdots \mathbf{B}(1){\bf v}_0.
\]
Then from \eqref{cascade}, for $n\ge 1$,
\[
\left[ \begin{array}{c} f(\bbb+{\aaa} n) \\ \vdots \\ f(\bbb+{\aaa} n+d-1)\end{array}\right] 
\ =\  \mathbf{A}(n)\mathbf{A}(n-1)\cdots \mathbf{A}(1){\bf v}.
\]
Since the hypotheses in the statement of Lemma \ref{lem: aal} are satisfied by $\mathbf{A}(\var)$ and
${\bf v}$, 
we also have that there exist rigid $p$-adic power series $f_1(\var),\ldots ,f_d(\var)$ such that,
for $n\ge 1$,
\[
\left[ \begin{array}{c} f_1(n) \\ \vdots \\ f_d(n)\end{array}\right] 
\ =\  \mathbf{A}(n)\mathbf{A}(n-1)\cdots \mathbf{A}(1){\bf v}.
\]
Thus for $n\ge 1$,
$$
f_1(n)\ =\ f(\bbb+{\aaa} n).
$$  
  By Strassman's theorem, 
$$
\{n\in \NN~:~f_1(n)=0\}
$$ 
is either finite or equal to $\NN$. Hence 
$$
\{n\in \NN~:~f(\bbb+\aaa n)=0\}
$$
 is either finite or equal to  $\NN$.  
 Let 
 $$
 X \ :=\ \big\{\bbb\in \{0,\ldots ,\aaa-1\}~:~f(\bbb+\aaa n) = 0 ~{\rm for~}n\in \mathbb{N}\big\}.
 $$  
 Then 
$$
\big\{n\in\mathbb{N}~:~f(n)=0\big\} \ =\  Z_0\cup \bigcup_{\bbb\in X} \big(\aaa \mathbb{N}+\bbb\big),
$$ 
where 
$Z_0$ is a finite set consisting of all $n$ such that $f(n)=0$ and $n\equiv \ccc~(\bmod \, \aaa)$ for some 
$\ccc\in \{0,\ldots ,\aaa-1\}\setminus X$.  The result follows.
\end{proof}

\section{Concluding remarks and examples}
We make some general remarks about Theorem \ref{thm: main} in this section.  

There are some difficulties that arise with trying to extend this argument to the collection 
of all holonomic sequences.  For example, if we take $$\displaystyle f(n)={2n\choose n},$$ then $f(n)$ is holonomic 
as it satisfies the polynomial linear recurrence
$$
(n+1) f(n+1) - 2(2n+1) f(n) \ = \ 0
$$ 
for $n\ge 0$.  Observe that there is no way to select a prime $p\ge 5$ and 
a natural number $a$ such that for all $\bbb\in\{0,\ldots ,\aaa-1\}$ we have 
a $p$-adic rigid analytic map $G_{\bbb}(\var)$ such that $G_{\bbb}(n)=f(\aaa n+\bbb)$ for all $n\in\mathbb{N}$.  
To see this, suppose that we have such a prime $p$ and a natural number $\aaa$ and let $\bbb=0$ and 
let $g(\var)$ be a $p$-adic rigid analytic map with $g(n)=f(\aaa n)$.  Since $g$ is continuous, 
there is some $j\ge 1$ such that $\big|g(\var+p^j)-g(\var)\big|_p<1/p$ for $\var\in \mathbb{Z}_p$.  
This then gives that
$$
f(\aaa p^jk)\ \equiv\  f(0)\ =\ 1~(\bmod\, p)
$$ 
for all natural numbers $k$.
But $f(pk)\equiv f(k)~(\bmod \, p)$ for all natural numbers $k$ and thus we must have
$$
f(\aaa k)\ \equiv\  1~(\bmod\, p)
$$ for all natural numbers $k$.
But if we choose $\ell$ such that $p^{\ell}\le \aaa <p^{\ell+1}$ and choose $k$ to be the 
smallest natural number such that $2\aaa k>p^{\ell+1}$ then we see that $p$ divides 
$f(\aaa k)$ using the formula
$$
\bigg|{2n \choose n}\bigg|_p \ = \ p^{\sum_{j\ge 1} 2\lfloor n/p^j\rfloor- \lfloor 2n/p^j\rfloor},
$$
and noting that each term appearing in the sum on the right-hand side is non-positive and 
that the term is strictly negative when $j=\ell+1$.  This is a contradiction.

In an earlier paper, the authors \cite{BBY} considered the set of zero coefficients to functions 
arising from solutions to a system of equations with certain prescribed properties.  In many cases,
 the solution sets were algebraic functions and hence differentiably finite.  
 The hypotheses guaranteed that the zero set was a finite union of arithmetic progressions 
 along with a finite set, although the methods used were very different and relied on studying 
 the behaviour of power series defined by systems of equations.  

It is also interesting to ask what happens in positive characteristic.  
Lech \cite{Le53} gave examples that showed that the conclusion of the Skolem-Mahler-Lech 
theorem does not hold if one eliminates the hypothesis that the field have characteristic $0$.  
For example if $p$ is a prime and $K=\mathbb{F}_p(t)$, 
then $f(n)=(1+t)^n-t^n-1$ satisfies a linear recurrence over $K$ but the zero set of $f$ is 
$\{1, p, p^2,\ldots \}$.  

Derksen \cite{Der} showed that the zero set $\mathcal{S}$ of a linear recurrence over a 
field of characteristic $p>0$ has the property that there is a finite-state automaton which takes 
as input the base $p$ expansion of a number $n$ and accepts the number if and only if 
$n\in \mathcal{S}$.  We call such sets $p$-\emph{automatic sets} 
(see the book of Allouche and Shallit \cite{AS} for a more precise definition).  
Recently, Adamczewski and the first-named author \cite{AB} have shown that the set of zero 
coefficients of an algebraic power series over a field of characteristic $p>0$ is $p$-automatic. 
 It thus seems reasonable to conjecture that if $f(n)$ is a holonomic sequence taking integer values 
 and $p$ is a prime number, then the zero set of the reduction of $f(n)$ (mod $p$) is a $p$-automatic set.  
 
 We also observe that $f(n)$ is a holonomic $\mathbb{Z}$-valued sequence that satisfies a polynomial-linear recurrence of the form 
 $$f(n) = \sum_{i=1}^d P_i(n)f(n-i),$$ then the argument we employed to prove that the zero set of $f(n)$ is a finite union of infinite arithmetic progressions along with a finite set can be proved under more general conditions.  In particular, it is sufficient that there exist infinitely many primes $p$ for which $P_d(x)$ does not have any roots modulo $p$.  The Chebotarev density theorem (see Lang \cite[Theorem 10, p. 169]{La94}) gives that this will occur precisely when there is some automorphism of the splitting field of $P_d(x)$ over $\mathbb{Q}$ whose natural action on the roots of $P_d(x)$ does not have any fixed points. 

Finally, we note that Theorem \ref{thm: main} is ineffective in the sense that 
if we know that our recurrence has only finitely many zeros, 
we cannot effectively bound the size of the largest zero in terms of data 
coming from the recurrence and the initial terms of the recurrence.  
Indeed, this is a notoriously difficult problem for linear recurrences and much work 
has been done on this problem by Evertse, Schlickewei, and Schmidt \cite{ESS}, 
who showed how one can obtain a quantitative version of the Skolem-Mahler-Lech theorem 
that bounds the number of exceptional zeros and the lengths of the arithmetic progressions 
in terms of such data.

\end{document}